\documentclass[12pt,leqno,twoside]{amsart}
\usepackage{amssymb}
\usepackage{latexsym}
\usepackage{verbatim}
\usepackage{epsfig}
\usepackage{changebar}

\newtheorem{theorem}{Theorem}[section]
\newtheorem{corollary}[theorem]{Corollary}
\newtheorem{lemma}[theorem]{Lemma}
\newtheorem{proposition}[theorem]{Proposition}
\theoremstyle{remark}
\newtheorem{remark}[theorem]{\sc Remark}
\theoremstyle{definition}
\newtheorem{definition}[theorem]{Definition}
\theoremstyle{remark}
\newtheorem{example}[theorem]{\sc Example}
\theoremstyle{remark}

\theoremstyle{remark}

\theoremstyle{remark}

\numberwithin{equation}{section}

\renewcommand{\Box}{_\square}    %\diamond

\newcommand{\cal}{\mathcal}

\newcommand{\h}{{\rm{ht}}}

\renewcommand{\int}{{\rm{int}}}
\newcommand{\gen}{{\rm{gen}}}

\newcommand{\Bif}{{\rm{Bif}}}

\newcommand{\Sing}{{\rm{Sing\hspace{2pt}}}}

\newcommand{\id}{{\rm{id}}}

\newcommand{\ity}{{\infty}}

\newcommand{\e}{\varepsilon}
\newcommand{\fisi}{{$\mathcal F${\rm{ISI}}}}
\newcommand{\fin}{\hspace*{\fill}$\Box$\vspace*{2mm}}

\newcommand{\cR}{{\cal R}}

\newcommand{\bC}{{\mathbb C}}

\newcommand{\bP}{{\mathbb P}}

\newcommand{\bZ}{{\mathbb Z}}
\newcommand{\bX}{{\mathbb X}}
\newcommand{\bY}{{\mathbb Y}}

\newcommand{\ft}{{\mathfrak{t}}}

\begin{document}

\title[Deformations of polynomials]
 {Deformations of polynomials, boundary singularities and monodromy}

\author{\sc Dirk Siersma}  % \ {\tiny and} \ Mihai Tib\u ar}

\address{D.S.: Mathematisch Instituut, Universiteit Utrecht, PO
Box 80010, \ 3508 TA Utrecht
 The Netherlands.}

\email{siersma@math.uu.nl}

\author{\sc Mihai Tib\u ar}

\address{M.T.:  Math\' ematiques, UMR-CNRS 8524,
Universit\'e de Lille 1, \  59655 Villeneuve d'Ascq, France.}

\email{tibar@agat.univ-lille1.fr}

%\thanks{}



%\date{\today}

\dedicatory{Dedicated to Vladimir Igorevich Arnol'd on the occasion of his 65th
anniversary}

%\commby{}

%%% ----------------------------------------------------------------------

\begin{abstract}
We study the topology of polynomial functions by deforming them generically. We explain how the non-conservation of the total ``quantity'' of singularity in the neighbourhood of infinity
is related to the variation of topology in certain families 
of boundary singularities along the hyperplane at infinity.
   
\end{abstract}

%%% ----------------------------------------------------------------------
\maketitle
%%% ----------------------------------------------------------------------

\setcounter{section}{0}
\section{Introduction}

 One of the fundamental tools for studying singularities 
 is deformation. In case of function germs with isolated singularity,  classification and deformation results have been found through the pioneering work of V.I. Arnol'd in the early 1970s, one of the most spectacular being the list of simple singularities, their possible deformations and their amazing connections to other branches of mathematics \cite{Ar1}. More recently, in the 1990s, the study of the behaviour of polynomial functions at infinity lead to a new branch: singularities at infinity of polynomials. In his 1998 paper \cite{Ar2}, Arnol'd gave a list of simple germs of fractions and polynomials at infinity, together with their specializations diagrams. One of the questions which may arise, placed in the global setting of polynomial functions, is the following: {\em how does the topology of the polynomial change when deforming it?} 
 
Since Broughton's fundamental paper \cite{Bro} there has been an enormous 
interest in studying the variation of topology of the fibers of a given polynomial, partly due to the relation to the Jacobian Conjecture.
In our paper, we study a polynomial by embedding it into a family. This gives a new viewpoint and non-trivial results occur even in cases where the topology at infinity of the original polynomial is constant (see \S \ref{s:arnold}).   
 
It is well-known that an isolated singular point of a holomorphic function can be deformed into finitely many singularities and that the sum of their Milnor numbers equals the original one. A generic deformation is a Morsification and it can be used to define vanishing cycles, see e.g. the Appendix of Brieskorn's paper \cite{Br}. 
 
 Unlike the case of germs, deformations of polynomials do not conserve  neither the total Milnor number, nor the number of vanishing cycles.   
      While focusing on the behaviour at infinity, there is a way to define generic deformations of $f_0$, within the class of degree $d$ polynomials, by asking that, for all $s\not= 0$, the total Milnor 
number $\mu(s)$ of $f_s$ is exactly equal to $(d-1)^n$.
Such polynomials $f_s$, which we call {\em general-at-infinity}, have also the maximum number of vanishing cycles (Proposition \ref{p:max}, Corollary \ref{l:bif}).
         
 Our study is placed in the framework of polynomials with isolated singularities, including singularities at infinity, in some appropriate sense.  We show that in this class, the number of vanishing cycles is upper semi-continuous (Proposition \ref{l:semicont}), which is just opposite as in the case of germs.
 
 The loss of vanishing cycles is actually due to the loss of some of the singular points of $f_s$, as $s\to 0$. 
This can happen in two ways: the corresponding critical value in $\bC$  of such a singularity tends either
  to a finite limit or to infinity. This phenomenon produces 
  singularities at infinity for $f_0$, which are quantified by the total Milnor-L\^e number $\lambda$. In other words, loss of $\mu$ yields gain of $\lambda$: 
Theorem \ref{t:lambdaexchange} shows how the  $\mu 
\rightleftarrows \lambda$ exchange takes place whenever singularities
tend to some point at infinity of a compactified fiber. In this way one can detect jumps in the topology by studying the behaviour of the affine critical points in the deformation. 
 
 Our approach is based on the interplay between the affine singularities of the polynomial and the singularities of its compactified fibers on the hyperplane at infinity $H^\ity$.  Actually there are finitely many one-parameter families of {\em boundary singularities} attached to the compactifies fibers of a polynomial function and their restrictions to the hyperplane at infinity. Arnol'd has already remarked in \cite{Ar2} that there is a close relation between the classification of meromorphic germs of type $H/Z^d$ and the one of germs of boundary singularities \cite{Ar1}.

       In a one-parameter deformation $f_s$ of $f_0$ there is a well-defined monodromy of a general fiber of $f_s$ (where $s$ is close enough to $0$) over a small circle $\partial \bar D \times \{ t\}$ centered at $(0,t)\in \bC \times \bC$. We call it the {\em $s$-monodromy}.  Depending on the position of $t$ with respect to the affine discriminant $\Delta_P$ (Definition \ref{d:gamma}), there are 3 types of $s$-monodromy, which we call {\em generic, atypical}, respectively {\em at infinity}.
       We compare these global monodromies to the monodromies of the germs of boundary singularities of $f_0$. On the other hand we decompose the 
    generic $s$-monodromy according to the types of branches of $\Delta_P$. We find in \S \ref{s:zeta} several formulas
    for the zeta-function of the monodromy.
    
    There is a whole section of examples on the $\mu 
\rightleftarrows \lambda$ exchange and the $s$-monodromy. Finally, we explain to what extent our study of a global situation (deformations of polynomials) relates to 
    Arnol'd's classification of simple germs of meromorphic functions of type $H/Z^d$. 
     
\hyphenation{For-schungs-ins-ti-tut}
{\sc Acknowledgements.}   This research originated in a discussion with Vladimir Igorevich Arnol'd in 1999 in Lille, on how vanishing cycles may appear in the list of simple meromorphic germs (see \S \ref{s:arnold}).  
Part of it was done as a 2002 RiP program at Oberwolfach, with support from the Volkswagen-Stiftung.
The authors benefited from the excellent working conditions
offered by the Mathematisches Forschungsinstitut Oberwolfach.    
      
%%%%%%%%%%%
%%%%%%%%%%%%%%%%%%%%%%

%%%%%%%%%%%%%%%
 
\section{Invariants associated to a polynomial functions}\label{s:inv}
%%%%%%%%%%%%%%%

For any polynomial function $f: \bC^n \to \bC$, there exists a minimal finite set 
$B_f\subset \bC$, called the set of {\em atypical values}, such that the restriction $f_|: 
\bC^n \setminus f^{-1}(B_f) \to \bC\setminus B_f$ is a locally trivial fibration. If $f$ 
has isolated singularities at infinity, in the sense of \cite{ST}, then its general fiber 
$G$ is homotopy equivalent to a bouquet of spheres $\vee_\gamma S^{n-1}$, cf. {\em loc. 
cit}. 

  In this case, the vanishing cycles
are quantified by two well-defined, non-negative integers \cite{ST}:
\[ \begin{array}{ll} \mu= \mbox{ the total Milnor number  of the affine singularities},\\
 \lambda= \mbox{ the total Milnor-L\^e number at infinity}
 \end{array} \]  
and one has:
%%%%%%%%%%%%%%
\begin{equation}\label{eq:mulambda}
 b_{n-1}(G) = \mu + \lambda.
\end{equation}

  In this paper we consider deformations of polynomials within the following class:
%%%%%%%%%%%%%%%%%%%%%%%%%%  
\begin{definition}\label{d:fisi}
 We say that $f$ has {\em fiberwise isolated singularities at infinity} (for short,  $f$ is 
a \fisi \  polynomial) if the projective closure of any fiber of $f$ and its slice by the hyperplane at infinity $H^\ity$ have at most isolated singularities.
\end{definition}
%%%%%%%%%%%%%%%%%%%%%%%%%%  
This class is large enough to include all reduced polynomial functions in 2 variables. 
It is a subclass of polynomials having isolated singularities at infinity, in the sense of \cite{ST}. 
So the numbers  $\mu$ and $\lambda$ are well-defined. 
The singularities at infinity of a \fisi \ polynomial are detected as follows. 

Let $\bX = \{ ([x;x_0],t)\in \bP^n\times \bC \mid \tilde f - tx^d_0 = 0\}$
 and let $\tau : \bX \to \bC$ be the projection to the $t$-coordinate, where $d= \deg f$ and $\tilde f$ is the degree $d$ homogenisation of $f$ with the variable $x_0$. Then the singular locus of $\bX$ consists of the lines $\{ p_i\} \times \bC$, where $p_i$ is a singular point of $\bX_t:= \tau^{-1}(t)$, situated on the part at infinity $\bX_t^{\ity}:= \tau^{-1}(t) \cap H^{\ity}$ (and not depending on $t$).

 Let now $(p_i,t)\in \bP^{n-1}\times \bC$ be a singular point of  $\bX_t^{\ity}$. This may  
 be a singular point of $\bX_t$, or a point where $\bX_t$ is non-singular but tangent to $H^{\ity}$ at $p_i$.
  Let then denote by:
  \[ \mu_{p_i}(t) := \mu(\bX_t, p_i), \ \ \ \mu_{p_i}^\ity(t) := \mu(\bX_t^{\ity}, p_i),\]
  the Milnor numbers of the respective hypersurface germs at $p_i$.
    Note that $\mu_{p_i}^\ity(t)$ is independent of $t$, therefore we shall use the 
notation $\mu_{p_i}^\ity$. In contrast, $\mu_{p_i}(t)$ may jump at a finite number of 
values of $t$. Let then denote by 
$\mu_{p_i, gen}$ the value of $\mu_{p_i}(t)$ for generic $t$.
 The jump of  $\mu_{p_i}$ at $t$ is measured by the local Milnor-L\^e number, cf. \cite{Pa}, \cite{ST}:
 %%%%%%%%%%%
\begin{equation}\label{eq:lambda}
\lambda_{p_i, t} = \mu_{p_i}(t) - \mu_{p_i, gen}.
\end{equation}
%%%%%%%%%%%%%%%%%%%
 We say that $(p_i,t)$ is a {\em singularity at infinity of $f$} if and only if $\lambda_{p_i, t}> 0$. Then the set of atypical values $B_f$ consists of the values $t\in \bC$ such that either $t$ is a critical value of $f$ or $(p_i,t)$ is a singularity at infinity of $f$ for some $p_i\in H^\ity$.
  The total Milnor-L\^e number of $f$ is the sum of the local ones:
\begin{equation}\label{eq:lambdatotal}
 \lambda = \sum_{t\in B_f}\sum_i \lambda_{p_i, t} =  \sum_{t\in B_f}\sum_i(\mu_{p_i}(t) - \mu_{p_i, gen}).   
\end{equation}
%%%%%%%%%%%%%%%%%%%
Let us compute the 
$(n-1)$-th betti number of the fibers $F_t = f^{-1}(t)$, using information on the singularities of their compactifications $\bX_t$ and of their restrictions to the hyperplane at infinity $\bX_t^{\ity}$. 
   Recall that the Euler characteristic of the general fiber of a generic polynomial of degree $d$ is $\chi(n,d) = 1+(-1)^{n-1}(d-1)^n$. 
  Comparing it to the Euler characteristic of $F_t$, one obtains the following formula, see e.g. \cite{Di}: 
\[
\chi(n,d) - \chi(F_t) = (-1)^{n-1}\sum_i(\mu_{p_i}(t) + \mu_{p_i}^\ity).
\]
 This is valid for all $t$, as long as $F_t$ is non-singular. When $F_t$ contains affine critical 
points, then the formula has to be corrected by adding to the right hand side the local 
Milnor numbers of these singularities. 
  When $F_t$ is the general fiber $G$, we obtain:
  
\begin{equation}\label{eq:chi}
  b_{n-1}(G) = (d-1)^n - \sum_i(\mu_{p_i, gen} + \mu_{p_i}^\ity).
 \end{equation}
 
We have now two different formulas, (\ref{eq:mulambda}) and (\ref{eq:chi}), both computing the number of vanishing cycles for a \fisi \ polynomial.
%
%%%%%%%%%%%%%%%%%%%%%%%%%%%%%
\section{Deformations of \fisi \ polynomials}
%%
%%%%%%%%%%%%%%%%%%%%%%%%%%%%%
%%%%%%%%%%%%%%%%%%%%%%%%%%%%%
In the following we consider a constant degree deformation of a polynomial 
 $f_0$, i.e. a family $P:\bC^n \times \bC^k \to \bC$  of polynomial functions  $P(x,s) = 
f_s(x)$ such that $\deg f_s =d$ for all values of the parameter 
$s$ in a small neighbourhood of $0\in\bC^k$.  We assume that $P$ is holomorphic or polynomial 
in the parameters $s$ and polynomial in the variables $x$.
 
%%%%%%%%%%%%
In this section we work with one-parameter deformations, i.e. $k=1$.
The following sets play an important role in our study:

\begin{equation}\label{eq:not}
\begin{array}{c}  W_s := \{ [x] \in \bP^{n-1} \mid   \frac{\partial P_d}{\partial 
x}(x,s)=0\},\\ \\
 \Sigma_s := \{ [x] \in \bP^{n-1} \mid   P_{d-1}(x,s)=0 \} \cap W_s,
 \end{array}
 \end{equation}
 where $P_q$ stays for the 
degree $q$ homogeneous part of $P$ in variables $x\in \bC^n$.

That $f_0$ is a \fisi \ polynomial is equivalent to: $\dim W_0 \le 0$ and $\dim \Sing f_0 
\le 0$. Remark that the second condition is not a consequence of the first one (example: $f_0 
= x^2$, as polynomial of $2$ variables).
%%%%%%%%%%%%%%%%%%%%%%%%%%
\begin{definition}\label{d:main}
We say that $P$ is a {\em \fisi \ deformation}
of $f_0$ if $P$ has the property that $\dim W_0 \le 0$ 
and that $\dim \Sing f_s \le 0$ for all $s$ in a small neighbourhood of $0$.
\end{definition}
%
%%%%%%%%%%%%%%%%%%%%%%%%%%%%%
Our definition of \fisi \ deformation
contains the condition that $f_0$ is \fisi \ and it implies
 that $f_s$ is \fisi \ for any $s$ in some small 
neighbourhood of $0$. 
%%%%%%%%%%%%%%%%%%%%%%
%%%%%%%%%%%%%%%%%%%%%%%%%%%%

\begin{example}
A constant degree deformation such that, for all $s$,  $f_s$ is a polynomial in 
2 variables with isolated singularities in the affine, is a \fisi \ deformation.
\end{example}

%%%%%%%%%%%%%%%%%%%%%%
\bigskip

We attach to the  family $P$ the following hypersurface:
\[ \bY = \{ ([x:x_0], s, t) \in \bP^n \times \bC \times \bC \mid \tilde 
P(x,x_0,s) - tx_0^d = 0 \},\] 
where $\tilde P$ denotes the homogenised of $P$ by the variable $x_0$, 
considering $s$ as parameter varying in a small neighbourhood of $0\in \bC$.
Let $\tau : \bY \to \bC$ be the projection to the $t$-coordinate. This 
extends the map $P$ to a proper one in the sense that $\bC^n\times \bC$ is embedded into $\bY$ and that $\tau_{|\bC^n\times \bC} = P$. Let $\sigma : \bY \to 
\bC$ denote the projection to the $s$-coordinates.
We shall use the notations 
$\bY_{s,\bC} := \bY\cap \sigma^{-1}(s)$, $\bY_{\bC, t} := \bY\cap 
\tau^{-1}(t)$ and 
$\bY_{s,t} := \bY_{s,\bC}\cap \tau^{-1}(t) = \bY_{\bC, t}\cap 
\sigma^{-1}(s)$. Note 
that $\bY_{s,t}$ is  
the closure in $\bP^n$ of the affine hypersurface $f_s^{-1}(t)\subset \bC^n$.
Let $\bY^\ity := \bY\cap \{ x_0=0\} = \{ 
P_d(x,s) = 0\}\times \bC$ be the hyperplane at infinity of $\bY$.

%%%%%%%%%%%%%%%%%

The singular locus of $\bY$, namely:
\begin{equation}\label{eq:singY}
 \Sing \bY = \{ x_0 =0, \frac{\partial P_d}{\partial x}(x,s) =0, \  P_{d-1}(x,s)=0, 
\frac{\partial P_d}{\partial s}(x,s)=0 \} \times \bC
\end{equation}
is a subset of $\bY^\ity$ and is a product-space by the $t$-coordinate. 
It depends only on the degrees $d$ and $d-1$ parts of $P$ with respect to the variables $x$.

%%%%%%%%%%%%%%%%%
For a fixed $s$ in a small neighbourhood of $0\in \bC$, the singular locus of $\bY_{s,\bC}$ is the  analytic set:

\begin{equation}\label{eq:sing}
\Sing \bY_{s,\bC} = \Sigma_s \times \bC \subset \bY_{s,\bC}\cap \bY^\ity.
\end{equation}

 It is the union, over all $t\in \bC$, of the singularities of 
the fibers $\bY_{s,t}$ situated on the part at infinity $\bY_{s,t}^\ity := \bY_{s,t} \cap H^\ity$. Besides these singularities, there might be points $(p,s,t)\in \bY^\ity$ where the fiber $\bY_{s,t}$ is not 
singular but it is tangent to the hyperplane at infinity $H^\ity$. Those are exactly the points such that
$p\in W_s\setminus \Sigma_s$. Notice that $W_s$ is the singular locus of 
$\bY^\ity_{s,t}$ (and does not depend on $t$). Let us set the following notations, which we shall frequently use:
%%%
\[\phi(x) := \frac{\partial P}{\partial s}(x,0), \ \  \phi_d := \mbox{ degree } d \mbox{ 
part of } \phi .
\]
%%%%%%%%%%%%%%%%%%%%%%
\begin{proposition}\label{l:phi}
Let $P$ be any deformation of a \fisi \ polynomial $f_0$. Then:

\begin{enumerate}
\item   $\bY$ is non-singular if and only if $\Sigma_0 \cap \{\phi_d =0\} = \emptyset$. In 
this case the critical set of the map $\sigma : \bY \to \bC$ is $\Sigma_0 \times \bC 
\subset\bY_{0,\bC}$.
\item  $\bY^\ity$ is non-singular if and only if $W_0 \cap \{\phi_d =0\} = \emptyset$. In 
this case the critical set of the map $\sigma_| : \bY^\ity \to \bC$ is $W_0 \times \bC 
\subset\bY^\ity_{0,\bC}$.
\end{enumerate}
\end{proposition}
\begin{proof}
(a) We have that $(\Sigma_0 \cap \{ \phi_d = 0\}) \times \bC = \Sing \bY \cap \{ \sigma 
=0\}$.
Since $\Sing \bY$ is a germ of an analytic set at $\bY_{0,\bC}$, it follows that $\Sing \bY$ 
is empty if and only if $(\Sing \bY) \cap \{ \sigma =0\}$ is empty. 
The critical locus of $\sigma : \bY \to \bC$ is a germ at $\bY_{0,\bC}$ and therefore it is 
$\Sing \bY_{0,\bC} = \Sigma_0 \times \bC$.
The proof of (b) is similar.
\end{proof}
%%%%%%%%%%%%%%%%%%%%%%%%%%%%%%
\begin{corollary}\label{c:phi}
Let $f_0$ be a \fisi \ polynomial and let $P$ be a deformation such that  $W_0 \cap \{ \phi_d = 0\}$ is 
empty. Then 
$\bY_{s,t}$ and $\bY^\ity_{s,t}$ are non-singular, for all $t\in \bC$ and all $s\not= 0$ 
close to $0$.
\end{corollary}
\begin{proof}
From Proposition \ref{l:phi}, it follows that $W_s = \emptyset$, for $s\not= 0$.
By (\ref{eq:sing}), we get that $\Sing \bY_{s,\bC}$ is empty and so is $\Sing \bY^\ity_{s,\bC}$. Since both spaces are products by the $t$-coordinate, we get our claim.
\end{proof}

%%%%%%%%%%%%%%%
Since we shall use this several times, we should emphasize that
in case $W_0 \cap \{ \phi_d = 0\}$ is empty, the spaces $\bY$ and $\bY^\ity$ are non-singular. This is in contrast to the compactifying spaces  $\bY_{0,\bC}$, resp. $\bY_{0,\bC}^\ity$, of the original polynomial 
$f_0$, which spaces have non-isolated singularities as soon as $\Sigma_0 \not= \emptyset$, resp. $W_0 \not= \emptyset$. For instance, in $\bY$
we have smoothings $\bY_{s,\bC}$ of $\bY_{0,\bC}$ and in $\bY_{\bC, t}$ we have smoothings $\bY_{s,t}$ of  $\bY_{0,t}$, for all
 $t$, for $s\not=0$.
%%%%%%%%%%%%%%%%%%%%%%%%%%
\begin{example}\label{ex:yomdin}
 {\em Yomdin deformations} of $f_0$ are of the form:
\[  \ P = f_0  - sx_1^d, \]
where $x_1$ denotes a coordinate. 
Sufficiently general Yomdin deformations, i.e. such that  $\Sigma_0 \cap \{ x_1 = 0\} = \emptyset$, are \fisi \ deformations. We can use 
the chart $x_1 =1$ to obtain a description of $\bY$ as graph of a function:
\[  s = (f_0)_d(1, x_2, \ldots, x_n) + x_0(f_0)_{<d} - t x_0^d,\]
where $(f_0)_{<d}$ denotes $(f_0)_{d-1} + x_0(f_0)_{d-2}+\cdots$. 
\end{example}
%%%%%%%%%%%%%%
%%%%%%%%%%%%%%%
%%%%%%%%%%%%%%
%%%%%%%%%%%%%%%

\subsection{Generic deformations}
%%%%%%%%%%%%%%
%%%%%%%%%%%%%% 
 \begin{definition}\label{d:gi}
A  polynomial $f : \bC^n \to \bC$ such that, for any $t\in \bC$, the compactified fiber 
$\overline{f^{-1}(t)} \subset \bP^n$ is transversal to the hyperplane at infinity $H^\ity = 
\bP^n\setminus \bC^n$, will be 
 called {\em general-at-infinity} (abbreviated GI).  We say that a 1-parameter deformation 
$P:\bC^n \times \bC \to \bC$ is GI if the polynomial $f_s$ is GI, for any $s\not= 0$ within 
a small neighbourhood of $0$.
 \end{definition}
We have the following property.
%%%%%%%%%%%%%%%%%%%%%%%%
\begin{proposition}\label{p:max}
 The general fiber $G$
of a general-at-infinity polynomial  $f$
is homotopy equivalent to a bouquet of $(n-1)$-spheres and the number of spheres is $\mu = 
b_{n-1}(G) = (d-1)^n$.
\end{proposition}
%%%%%%%%%%%
\begin{proof} 
Since $f$ is GI,  all its fibers are transversal
to the hyperplane at infinity. This implies that 
$f$ has only isolated singularities in the affine and 
that the Milnor numbers $\mu_{p_i}(t)$ and  $\mu_{p_i}^\ity$ are all zero, for all $t$. Then formula (\ref{eq:chi}) shows that $b_{n-1}(G) = (d-1)^n$.
Further, since our $f$ is a \fisi\ polynomial, it is in particular a polynomial with at most isolated singularities at infinity, in the sense of \cite{ST}, and therefore its general fiber is a bouquet of $(n-1)$-spheres,
by \cite[Theorem 3.1]{ST}.
\end{proof}
%%%%%%%%%%%%%%%%%%%
\begin{corollary}\label{l:bif}
A \fisi \ polynomial $f : \bC^n \to \bC$ is general-at-infinity if and only if $b_{n-1}(G) = (d-1)^n$. 
\end{corollary}
\begin{proof}
In one direction, this follows from Proposition \ref{p:max}. The reciprocal is a direct 
consequence of the formula (\ref{eq:chi}) for $b_{n-1}(G)$ in terms of $\mu_{p_i}$ and 
$\mu_{p_i}^\ity$. 
\end{proof}
%%%%%%%%%%%%%
%%%%%%%%%%%%%%
  A GI polynomial may have different configurations of isolated 
singularities in the 
affine, provided that their total Milnor number is precisely $(d-1)^n$.
A large class of GI deformations is given by the following result.
%%%%%%%%%%%%%%%
\begin{corollary}\label{c:nongi}
Let $P$ be a deformation of a \fisi \ polynomial $f_0$ such that $W_0 \cap \{ \phi_d = 0\} 
=\emptyset$. Then $P$ is general-at-infinity.
\end{corollary}
\begin{proof}
 By definition, $P$ is GI 
  if and only if  $\bY_{s,t}$ and $\bY^\ity_{s,t}$ are non-singular, for all $t\in \bC$ and 
all $s\not= 0$ close to $0$.
 Then our claim follows from Corollary \ref{c:phi}.
\end{proof}
%%%%%%%%%%%%%%%%%%%%%%%%%%%
%%%%%%%%%%%%%%%%%%%%%%%%%%%%%%%%%%%%%%%%%%%%

%%%%%%%%%%%%%%%
%%%%%%%%%%%%%%%%%%%%%%
%%%%%%%%%%%%%%%%%%%%
%%%%%%%%%%%%%%%%%%

\section{The affine critical locus}
%%%%%%%%%%%%%%%%%%%%
%%%%%%%%%%%%%%%%%%%%
In case of several parameters, we have of course more room for 
constructing generic deformations. For instance, let $P\colon \bC^n\times \bC^n \to \bC$ be a deformation of $f$ such that $P_d= f_d + 
s_1x_1^d + \cdots + s_nx_n^d$. Then $\bY$ is nonsingular and 
$\bY_{\bC,t}$ is nonsingular too, for any $t$.

Let us now focus on the critical locus of the deformation. 

\begin{definition}\label{d:gamma}
Let $\Psi := (\sigma, \tau) :  \bY\setminus \bY^\ity \to \bC^k \times \bC$ be the map induced by the couple of
projections. We define the {\em affine critical locus} $\Gamma_P$  to be the 
 closure in $\bP^n \times \bC^k \times \bP$ of the set of points where  $\Psi$ is not a 
submersion. Let us call {\em affine discriminant} and denote by $\Delta_P$ the closure in 
$\bC^k \times \bP$ of the image $\Psi( \Gamma_P \cap (\bC^n \times \bC^k \times \bC))$ of the affine part of $\Gamma_P$.
\end{definition}
%%%%%%%%%%%%%%%%%%%%%%%
Since we are only interested in the germ of $\Gamma_P$ at $\bP^n \times \{ 0\}\times \bP$, 
respectively in the germ of $\Delta_P$ at $\{ 0\}\times \bP$, we shall tacitly use the same 
notations for the germs. 
 
It follows from the definition that $\Gamma_P$ is a closed analytic set and its affine part $\Gamma_P  \cap 
(\bC^n\times \bC^k \times \bC)$ is just the union, over $s\in \bC^k$, of the affine critical loci of 
the polynomials $f_s$. 
 Note that, if 
$\dim \Sing f_s \le 0$ for any $s$, then $\dim \Gamma_P \le k$.
%%%%%%%%%%%%%%%%%%%%%%%%%%%%%%%%%%
\subsection{Affine critical locus of \fisi \ deformations}
%%%%%%%%%%%%%%%%%%%%%%%%%%%%%%%%%%%
%%%%%%%%%%%%%%%%%%%%%%%%%%%%%%%%%
 Let now assume that $P$ is a one-parameter \fisi \  deformation. It follows  that, if 
$\Gamma_P$ is not empty, then it is a curve.

An irreducible component $\Gamma_i$ of the germ $\Gamma_P$ can be of one of the following three types:
\begin{itemize} 
\item[-]Type I: $\Gamma_i \cap \{\sigma =0\} \in \bC^n \times \{ 0\} \times  \bC$. Then 
$\Gamma_i$ is not in the neighbourhood of infinity.
\item[-] Type II:  $\Gamma_i \cap \{\sigma =0\} \in \bY^\ity$, a point having finite 
$t$-coordinate. 
\item[-] Type III:  $\Gamma_i \cap \{\sigma =0\}$ is a point having infinite $t$-coordinate. 
\end{itemize}
 Let us decompose $\Gamma_P = \Gamma_I \cup \Gamma_{II} \cup \Gamma_{III}$ according to 
these types.
%%%%%%%%%%%%%%%%%%%%%%%%%%%%%%%%%
If the total Milnor number $\mu(s)$ is not constant in the deformation $P$, then necessarily
$\Gamma_{II} \cup \Gamma_{III} \not= \emptyset$. The points   where $\Gamma_{II}$ 
intersects the hyperplane at infinity are shown to play
an important role. The similar question for $\Gamma_{III}$ requires the completion with 
$t=\ity$. 
%%%%%%%%%%%
\begin{lemma}\label{l:curve2}
Let $P$ be a one-parameter \fisi \  deformation. Then $\Gamma_{II}\cap  \{ x_0= 0 \}\subset 
\Sigma_0 \times \{ 0\} \times\bC$ and $\Gamma_{III}\cap  \{ x_0= 0 \}\subset W_0\times \{ 0\} \times \{ \ity\}$. 
\end{lemma}
\begin{proof}
From the definition of $\Gamma_P$ it follows that $\Gamma_P \cap \{ \sigma =s\} \subset 
\bC^n \times \bC$, for any $s\not= 0$. Therefore $\Gamma_P \cap \{ x_0= 0 \}\subset H^\ity 
\times \{ 0\} \times \bC$.
The curve $\Gamma_P \cap \bY$ is included into the critical set of the map $\bar \Psi = 
(\sigma, \tau) : \bY \to \bC\times \bP$.
 Slicing this with $\{ x_0= 0 \}$ and $\{ \sigma =0\}$ gives $\{ \frac{\partial 
P_d}{\partial x} (x,0) = 0, \ P_{d-1}(x,0) = 0\} \times \bC$.
 
 For $\Gamma_{III}$, we homogenise the equations of $\Gamma_P \cap (\bC^n \times \bC \times 
\bC)$, which are 
$\{ \frac{\partial P}{\partial x} (x,s) = 0, \ P(x,s) = t\}$, and get an algebraic set into 
$\bP^n \times \bC \times \bP^1$, which we slice by $t=\ity \in \bP^1$ and by $\{ \sigma =0\}$. This gives the set 
$\{ \frac{\partial f_d}{\partial x} (x) = 0\}\subset H^\ity$.
\end{proof}
%%%%%%%%%%%
 For general-at-infinity deformations we can prove more precise results, as follows (see 
also Theorem \ref{t:lambdaexchange}). 
\begin{proposition}\label{l:iomdin}
Let $P$ be a general-at-infinity deformation of a \fisi \ $f_0$.
If for some $t\in \bC$, $\bY_{0,t}$ is non-singular but tangent to $H^\ity$ then  $\Gamma_{III} \not= \emptyset$.
\end{proposition}
\begin{proof}
Our hypothesis on $\bY_{0,t}$ is equivalent to $\Sigma_0 =\emptyset$ and $W_0 \not= 
\emptyset$. Since $\Gamma_P \cap \bY_{0,t} \subset \Sigma_0 \times \bC$, it follows that 
there are no type II components of $\Gamma_P$. Moreover, $\Sigma_0 =\emptyset$ implies that 
$\lambda(0) = 0$. Since $W_0 \not= \emptyset$, the total Milnor number $\mu(0)$ is 
strictly less than $\mu(s)$, for $s\not= 0$, by (\ref{eq:chi}) and (\ref{eq:gammas}). 
This can only be due to the fact that affine critical values of $f_s$ tend to infinity, as 
$s\to 0$. This means that branches of type III do exist. 
\end{proof}
%%%%%%%%%%%%%%%
 %%%%%%%%%%%%%%%%%%%%%%
\begin{proposition}\label{l:tang}
Let $P$ be a one-parametre \fisi \ deformation and assume that the total space $\bY$ is non-singular.
If for some $t\in \bC$, $\bY_{0,t}$ has singularities on $H^\ity$ then all branches of 
$\Delta_P$ of type II are tangent to the $\{\sigma = 0\}$-axis.
\end{proposition}

\begin{proof}
Our condition on $\bY_{0,t}$ is equivalent to $\Sigma_0 \not=\emptyset$ and is a necessary 
condition for the existence of type II branches of $\Gamma_P$, see Lemma \ref{l:curve2}. 
Let $v= (p,t_0)\in \Gamma_P \cap \bY^\ity$. Denote by $\Gamma_{P,v}$ the germ of $\Gamma_P$ at $v$. 
We may assume that $\Gamma_{P,v}$ is an irreducible germ, if not, then we take an irreducible 
component of it. Take a parametrization of $\Gamma_{P,v}$, $\gamma(\e)$, with $\gamma(0)=v$.
 One should prove that the fraction $\sigma (\gamma(\e))/ \tau(\gamma(\e))$ tends to 0 as 
$\e\to 0$. This is indeed so since 
 $\tau$ is a linear function germ of the non-singular space $\bY$ and $\sigma$ has a 
singularity at $v$, hence it is in the square of the maximal ideal of holomorphic germ at 
$v$. Then our claim follows 
 by L\^e's result \cite[Proposition 1.2]{Le-c}.
\end{proof}
%%%%%%%%%%%%%%%

%
%%%%%%%%%%%%%%%
%%%%%%%%%%%%%%%%%%%%%%%%%%%%%% 

\section{Changes of topology within deformations}
%
%%%%%%%%%%%%%%%%
%%%%%%%%%%%%%%%
We study the change of topology within the class of \fisi \ polynomials.
Let $P$ be a \fisi \ deformation of a polynomial $f=f_0$ of degree $d$. We recall the 
notations in \S \ref{s:inv} and adapt them to our family $f_s$, as follows:
$G_s$ is the general fiber of $f_s$, $F_{s,t}:= f_s^{-1}(t)$,  $\mu(s)$ is the total Milnor 
number and $\lambda(s)$ is the total Milnor-L\^e number at infinity of $f_s$.
Note that if $f_s$ is GI then $\lambda(s) =0$ and $\chi(F_{s,t}) = \chi(n,d)$ if $t\not\in B_{f_s}$. We first have the following semi-continuity result: 
%%%%%%%%%%%%%%%%%%%%%%%%%%%%%%%%%%%
\begin{proposition}\label{l:semicont}
Let $P$ be a \fisi \ deformation of $f_0$. Then the number of 
vanishing cycles $\gamma(s) := \mu(s) +\lambda(s)$ is upper semi-continuous, i.e.  
\[ \mu(s) +\lambda(s) \ge \mu(0) +\lambda(0).\]
\end{proposition} 
%%%%%%%%%%%%%%%%%%%%%
\begin{proof}
 When examining the formula (\ref{eq:chi}) applied to $f_s$, namely:
 \begin{equation}\label{eq:gammas}
   \gamma(s) = b_{n-1}(G_s) = (d-1)^n - \sum_{p_i\in W_s}(\mu_{p_i, gen}(s) + \mu_{p_i}^\ity(s)),
\end{equation}  
we see that $\mu_{p_i, gen}(s)$ and $\mu_{p_i}^\ity(s)$ are lower semi-continuous since
 local Milnor numbers of hypersurface singularities.
\end{proof}
%%%%%%%%%%% 
From (\ref{eq:gammas}) and the semi-continuity of its $\mu$-ingredients, we also get: 
%%%%%%%%%%%%%%%%%%%%%%
\begin{corollary}\label{c:semi}
If $P$ is a $\gamma(s)$-constant deformation, then both $\mu_{p_i, gen}(s)$ and 
$\mu_{p_i}^\ity(s)$ are constant for small $|s|$.
\fin
\end{corollary}
%%%%%%%%%%%

Let us denote by $\mu_I(s)$, $\mu_{II}(s)$, $\mu_{III}(s)$ the sum of the Milnor numbers of 
the singularities of $f_s$ at the affine points of $\Gamma_I$,  $\Gamma_{II}$, $\Gamma_{III}$ respectively. We have $\mu(s) = \mu_I(s) + \mu_{II}(s) + \mu_{III}(s)$ 
and, if the deformation is general-at-infinity, then we know that $\mu(s) = (d-1)^n$.
%%%%%%%%%%%%%%%%%%%%
\begin{lemma}\label{c:mu23}
If $P$ is a general-at-infinity deformation then 
\[ \mu_{II}(s) + \mu_{III}(s) = \lambda(0) + \sum_{p_i \in W_0}\beta_{p_i}(0),\]
where $\beta_{p_i}(s) := \mu_{p_i,gen}(s) + \mu_{p_i}^\ity(s)$.
\end{lemma}
\begin{proof}
Since $f_s$ is GI for $s\not= 0$, we have $\lambda(s)=0$ and 
$W_s =\emptyset$, hence $\beta_{p_i}(s) = 0$. From the formula (\ref{eq:mulambda}) for 
$f_0$ and $f_s$, we get $b_{n-1}(G_s) - b_{n-1}(G_0) = \mu_{II}(s) + \mu_{III}(s) - \lambda(0)$, whereas from  
(\ref{eq:chi}), see also (\ref{eq:gammas}), we get $b_{n-1}(G_s) - b_{n-1}(G_0) = \Delta \beta := 
\sum_{p_i \in W_0}\beta_{p_i}(0)$.  
\end{proof}
%%%%%%%%%%%%

  Related to the semi-continuity result Proposition \ref{l:semicont}, a question which arises next would be: what is the behaviour of 
  $\mu(s)$ and of $\lambda(s)$ when $s\to \ity$?
      It is clear that the total Milnor number $\mu(s)$ cannot increase.  It can either be 
constant (e.g. in homogeneous deformations of polynomials)
or decrease, like in Example \ref{e:broughton}. 
In case $\mu(s)$ decreases, as $s\to 0$, we say that there is {\em loss of $\mu$ at infinity}; 
this happens if and only if $\Gamma_{II}\cup \Gamma_{III} \not= 0$.

In contrast, $\lambda(s)$ may either increase or decrease (or be constant). Examples where it 
increases may be found in Siersma-Smeltink paper \cite{SS}, see also Example \ref{e:mixed}. Examples of $\lambda$-constant deformations are all those with $\lambda(s)=0$.
  Example \ref{e:decrease} (of a non-\fisi\ deformation) shows that $\lambda(s)$ may decrease.
%%%%%%%%%%%%%
%%%%%%%%%%%%%
    An important fact is that we can now detect $\lambda$ by $\Gamma_P$.
%%%%%%%%%%%%%%%%
\begin{theorem}\label{t:lambdaexchange}
Let $P$ be a one-parameter \fisi \ deformation of a polynomial $f_0$ such that the total space $\bY$ is non-singular.  
Let $p\in \Sigma_0$.
 Then $(p,0,t_0) \in \Gamma_P$ if and only if $\lambda_{p,t_0} (0) \not= 0$ and in this case we have:
\[ \lambda_{p,t_0}(0) < \sum_{(a,s,t) \in \Gamma_P \cap \bY_{s,\bC}} \mu (\bY_{s,t}, a),\]
where $s\in D^*$ is fixed, close enough to $0$, and the sum is taken over the singular points of $f_s$ which tend to $(p,0,t_0)$ as $s\to 0$.
\end{theorem}
%%%%
\begin{proof}  
Consider the germ of $\bY$ at $(p,0,t_0)$ and the restriction map $\bar\Psi= (\sigma, \tau) : 
(\bY, (p,0,t_0)) \to \bC^2$.  The condition $(p,0,t_0) \in \Gamma_P$ implies that $\Sigma_0 
\not= \emptyset$, which implies that  $\sigma : \bY \to \bC$ has non-isolated singularities 
at $(p,0,t_0)$. 
Nevertheless $\tau$ is a linear function on $\bY$ and has isolated singularities on 
$\bY_{0,\bC}$.
Then the Milnor fiber of $\tau$ is contractible. After \cite{Le-m}, we may take a small ball $B$ in $\bY$ at 
$(p,0,t_0)$ and a 2-disk $D\times D'$ in $\bC^2$ at $(0,t_0)$, such that
 $\bar\Psi$ induces a locally trivial fibration $\bar\Psi_| : (B\cap \bY) \setminus (\Gamma_P \cup 
\{\sigma =0\} ) \to (D\times D') \setminus \Delta_P$.
Let $Y_{U,V} := B\cap \bY_{U,V}$, for $U\subset D$ and $V\subset D'$. Since $\tau$ is a linear function, we have that $Y_{D,t}$ is 
contractible, for any $t\in D'$.
Taking into account the singularities of the fibers of the projection $\sigma_| : Y_{D,t} \to D$ and using L\^e's attaching results \cite{Le-c}, we deduce the following homotopy equivalences:
%%%%%%%%%%%%%
\[ Y_{D,t_0} \stackrel{\h}{\simeq} Y_{s,t_0} \bigcup \mu_p(t_0) \mbox{ cells of dimension } n-1,\]
\[ Y_{D,t} \stackrel{\h}{\simeq} Y_{s,t} \bigcup \{\mu_p(t) \mbox{ cells of dimension } n-1 \}\bigcup 
\{ \alpha \mbox{ cells of dimension } n-1 \},\]
%%%%%%%%%%%%%
where 
\[ \mu_p(t) := \mu(\bY_{0,t},p), \hspace{1cm} \mu_p(t_0):= \mu(\bY_{0,t_0},p),\]
\[ \alpha := \sum_{(a,s,t) \in \Gamma_P \cap Y_{D,t}} \mu (\bY_{s,t}, a).\]

\begin{figure}[hbtp]
\begin{center}
\epsfxsize=5cm
\leavevmode
\epsffile{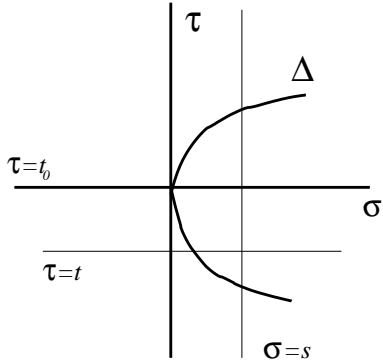}
\end{center}
\caption{{\em 
 Attaching cells}}
\label{f:0}   
\end{figure}

Since $Y_{s,t_0}\stackrel{\h}{\simeq} Y_{s,t}$ we have the equality 
$\mu_p(t_0)=  \mu_p(t) + \alpha$ and therefore
 $\alpha = \mu_p(t_0)- \mu_p(t) = \lambda_{p,t_0}(0)$, by (\ref{eq:lambda}).  
 Since $\Gamma_P$ consists  
of trajectories of affine singular points,  $\alpha$ is an intersection number, namely:
\[ \alpha = \int_{(p,t_0)}(\Gamma_P, \{ \tau = t_0\} ) = \int_{(0,t_0)}(\Delta_P, \{ \tau 
= t_0\} ),\]
where in the first intersection number, each point is counted with multiplicity equal to the Milnor number of 
the singularity of the fiber of $\bar \Psi$ at that point. The second intersection number involves 
in addition the multiplication by the degree of the covering $\Gamma \stackrel{\bar\Psi}{\to} 
\Delta$, for each irreducible component of $\Gamma$ and of $\Delta$.
This shows that $\alpha > 0$, which is our first claim. Let us emphasize that $\lambda_{p,t_0}(0)$ was identified with the intersection number between $\Gamma_{II}$
and the hyperplane  $\{ \tau = t_0\}$.

To prove the second part of the statement, let us first observe that the right-hand side of the claimed inequality is equal to the following intersection numbers, taken with appropriate multiplicities:
\[ \int_{(p,t_0)}(\Gamma_P, \{ \sigma = 0\} ) = \int_{(0,t_0)}(\Delta_P, \{ \sigma = 0\} ) 
.\]
By Proposition \ref{l:tang}, $\Delta_P$ is tangent to $\{ \sigma = 0\}$. This tangency implies 
that:
\[ \int_{(0,t_0)}(\Delta_P, \{ \sigma = 0\} ) > \int_{(0,t_0)}(\Delta_P, \{ \tau = t_0\} 
),\]
which gives the claimed inequality.
\end{proof}
%%%%%%%%%%%%%%%%%%%%%%%%%%%%
\begin{corollary}\label{c:a1}
In the conditions of Theorem \ref{t:lambdaexchange}, if 
$\deg f_0 > 2$ and the singularities of some fiber $\bY_{0,t}$ at $H^\ity$
are of type $A_1$ only, then $\Gamma_{II} = \emptyset$. 
\end{corollary}
\begin{proof}
Since $\deg f > 2$, the type $A_1$ of the singularities at $H^\ity$ of $\bY_{0,t}$ is 
independent of $t$. Consequently, there is no jump of $\mu$ at infinity. It follows, by formulas 
(\ref{eq:lambda}) and (\ref{eq:lambdatotal}), that $\lambda(0)=0$.
By Theorem \ref{t:lambdaexchange}, this means that there is no $\Gamma_{II}$.
\end{proof}

%%%%%%%%%%%%%%%%%%%%%%%%%%%%%%%%%%%%
%%%%%%%%%%%%%%%%%%%%%%%%%%%%%%%%%%%%%%%%%%%%%%%%%%%%%%%
%%%%%%%%%%%%%%%%%%%%%%%%%%%%%%%%%%%%%%%%%%%%%%%%%%%%%%%
%%%%%%%%%%%%%%%%%%%%%%%%%%%%%%%%%%%%%%%%%%%%%%%%%%%%%%%
%%%%%%%%%%%%%%%%%%%%%%%%%%%%%%%%%%%%%%%%%%%%%%%%%%%%%%%
%%%%%%%%%%%%%%%%%%%

\section{The $s$-monodromy and boundary singularities}\label{s:mono}

%%%%%%%%%%%%%%%%%%%%%%%%%%%%%%%%%%%%

%%%%%%%%%%%%%%%%%%%%%%%%
%%%%%%%%%%%%%%%%%%%%%%%%
%%%%%%%%%%%%%%%%%%%
Let $P$ denote a one-parameter, general-at-infinity deformation of a (non-GI) polynomial 
$f_0$ and define the {\em topological bifurcation diagram} in the $(s,t)$-target space of the map $\Psi=(\sigma, 
\tau) : \bY\setminus \bY^\ity \to \bC \times \bC$, as follows:
%%%%%%%%%%%%%%%%%
\begin{equation}\label{eq:topbif}  
 \Bif^{top}_P := \{ (s,t) \mid \Psi \mbox{ is not a topologically trivial fibration at } (s,t) \}\subset \bC \times \bC.
\end{equation}
%%%%%%%%%%%%%%%%%%%%%%%%%
  Since $P$ is a GI deformation, $\Bif^{top}_P$ is the union of the $t$-axis $\{ 
\sigma =0\}$ and the affine discriminant $\Delta_P \cap (\bC \times \bC)$.  

Let $F_{A,B}$ denote the set $(\bC^n\times \bC )\cap \sigma^{-1}(A)\cap 
\tau^{-1}(B)$, for some subsets $A\subset \bC$ and $B\subset \bC$.
 Take a small disc $D$ centered at the origin in the $s$-coordinate space 
$\bC$. The {\em $s$-monodromy} is the monodromy of the fibration over the 
circle $F_{\partial \bar D, \bC} \stackrel{\sigma}{\to}\partial \bar D$.
 
 We study this monodromy action on the general fiber $F_{s,t}$, where 
$(s,t)\not\in \Bif^{top}_P$, in relations to the singularities at infinity of 
$f_0$.
 %%%%%%%%%%%%%%%%%%% 
 
\begin{figure}[hbtp]
\begin{center}
\epsfxsize=3.5cm
\leavevmode
\epsffile{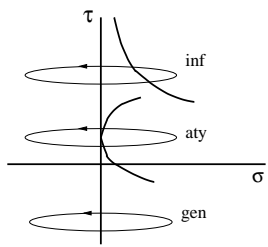}
\end{center}
\caption{{\em 
 Types of monodromy}}
\label{f:00}   
\end{figure}

 Let then $s\in\partial \bar 
D$ and let $h_{s,t}$ denote
 the $s$-monodromy over $\partial \bar D$ acting on the homology $H_*(F_{s,t},  \bZ)$. For a fixed value of $t$, there are 3 types of $s$-monodromies, corresponding to the 3 types of loops (see Figure \ref{f:00}): ``gen'' -- for a generic value of $t$, 
``aty'' -- for an atypical value of $t$, 
``inf'' -- for $t$ near infinity, in case of branch of type III.

In the first two cases we avoid intersections of $\bar D^*\times \{t\}$  with 
$\Delta_P$, in the third 
case we require intersections with branches of type III.

%%%%%%%%%%%%%%%%%%%%%%%%%%%
%%%%%%%%%%%%%%%%
%%%%%%%%%%%%%%%%
\subsection{The generic case}\label{ss:genmono}
%%%%%%%%%%%%%%%%
%%%%%%%%%%%%%%%%
Let $t\in \bC$ and let $D$ be a small enough disc as above, centered at $0\in \bC$, such that the disc 
$\bar D\times \{t\}$ does not intersect $\Delta_P$. 
The fiber $F_{s,t}$ has a geometric $s$-monodromy over $\partial \bar D$ which can be extended 
to a diffeomorphism of $F_{D,t}$ into itself, the latter being isotopic to the identity.
  The induced algebraic monodromy is called {\em generic 
$s$-monodromy} and will be denoted by $h_\gen$.
  
 Since $f_0$ is \fisi \ and the deformation $P$ is general-at-infinity,
 the fibers of the restriction to $D\times \{ t\}$ of the map $\bar\Psi : \bY \to \bC \times \bC$ have no singularities at $H^\ity$, 
except for the fiber $\bY_{0,t}$. So let $p_i\in \bY_{0,t}^\ity$
 be such a singularity and let $B_i \subset \bY$ be a Milnor ball at $p_i$.  
We assume that $D$ is small enough such that $(B_i, D)$ is Milnor data for any 
$i$.
 Then, by excision, we have the following direct sum decomposition of the relative homology:
\begin{equation}\label{eq:sum}
    H_*(F_{D,t}, F_{s,t}) =  \oplus_i H_*(B_i\cap F_{D,t}, B_i\cap F_{s,t}). 
\end{equation} 
 %%%%%%%%%%%%%%%%%%%%
\begin{lemma}\label{l:ball}
Let $P$ be a one-parameter \fisi \ deformation of a polynomial $f_0$ such that $W_0 \cap \{ \phi_d =0\} = \emptyset$. Then, for $n \ge2$:
\begin{enumerate}
\item $H_*(F_{D,t}, F_{s,t})$ is concentrated in dimension 
$n$.
\item
 $H_q(B_i\cap F_{s,t})$ is trivial if $q \not= 0, 1, n-1$, it is $\bZ$ if $q = 0$ and if $q= 1\not= n-1$, and it is isomorphic to $H_n(B_i\cap F_{D,t}, B_i\cap F_{s,t})$ if $q=n-1 \not= 1$. If $n=2$ then $H_1(B_i\cap F_{s,t}) \simeq \bZ \oplus H_2(B_i\cap F_{D,t}, B_i\cap F_{s,t})$.
\end{enumerate}
\end{lemma}
\begin{proof}
(a).   Since the family of hypersurfaces $\bY_{s,t}$, for $s\in D$ is a 
smoothing of the hypersurface with only isolated singularities $\bY_{0,t}$, it 
follows from \cite{Ti}, using \cite{HL}, that the relative homology 
$H_*(B_i\cap F_{D,t}, B_i\cap F_{s,t})$ is concentrated in dimension $n$. Our 
claim then follows by the direct sum decomposition (\ref{l:ball}). 
Note that the point (a) only requires the \fisi \ condition.

\noindent
(b).  By Proposition \ref{l:phi}, our space $\bY_{D,t}$ is nonsingular.
    Hence $B_i\cap \bY_{D,t}\setminus \bY^\ity$ is homeomorphic to a ball of dimension $n$ from which one takes out a
   coordinate hyperplane, so this space is homeomorphic to $\bC^n \setminus 
\bC^{n-1}$. Then (b) follows from the long exact sequence of the pair 
$(B_i\cap F_{D,t}, B_i\cap F_{s,t})$.
 \end{proof}  
%%%%%%%%%%

Let $h_i$ denote the $s$-monodromy action on  $H_{n-1}(B_i\cap F_{s,t})$. This 
can be interpreted as the Milnor monodromy of a {\em boundary singularity}. Indeed, 
the function $\tau$ on the ball $B_i\cap \bY_{D,t}$ has an isolated 
singularity at $p_i\in \bY^\ity_{0,t}$. Its restriction $\tau_{|\bY^\ity}$ to the non-singular hyperplane 
$\bY^\ity$ has also an isolated singularity at $p_i$. It was Arnol'd who 
founded the study of germs of boundary singularities and classified simple germs, see \cite{AGV}. In his paper \cite{Ar2}, Arnol'd shows that the classification of meromorphic germs of type $H/Z^d$ is closely related to the classification of boundary singularities (see also \S \ref{s:arnold}). 

To the pair of function germs $(\tau, \tau_{|\bY^\ity})$ at $p_i$, one associates the 
pair of Milnor fibers $(B_i\cap \bY_{s,t}, B_i\cap \bY^\ity_{s,t})$.
In a subsequent work, A. Szpirglas \cite{Sz} studied 
vanishing cycles of boundary singularities, from the point of view of a certain 
duality between the pair of Milnor fibers and their difference. After \cite{Sz}, there are 
two Lefschetz-dual short exact 
sequences, which in our case take the following form (in case $n\ge 2$):
%%%%%%%%%%%%%%%
\begin{equation}\label{eq:dual}
\begin{array}{ccccccccc}
0 & \to & H_{n-1}(B_i\cap \bY_{s,t}) & \to & H_{n-1}(B_i\cap \bY_{s,t}, 
B_i\cap \bY^\ity_{s,t}) & \stackrel{\partial}{\to} & \tilde H_{n-2}(B_i\cap 
\bY^\ity_{s,t}) & \to & 0 \\ \\
0 & \leftarrow & H_{n-1}(B_i\cap \bY_{s,t}) & \leftarrow  & H_{n-1}(B_i\cap 
F_{s,t}) & \stackrel{L}{\leftarrow} & \tilde H_{n-2}(B_i\cap \bY^\ity_{s,t}) & 
\leftarrow  & 0, 
 \end{array}
\end{equation} 
where $L$ is the Leray ``tube'' map.
%% Argument for duality: similar as in Broughton.
%%%%%%%%%%%%%%%%%%
%%%%%%%%%%%%%%%%%%

\section{Zeta-function of the $s$-monodromy}\label{s:zeta}
%%%%%%%%%%%%%%%%%%%
%%%%%%%%%%%

 For a monodromy action $h$ on a space $M$ one defines the {\em zeta-function} on  
 the $\bZ$-homology $H_*(M)$, as follows:
 \[ \zeta_M(\ft) := \prod_{j\ge 0} [ \det (I-\ft\cdot h \mid 
H_j(M))]^{(-1)^{j+1}}.\]
  With this sign convention, the degree of the zeta-function is equal to $-\chi (M)$.
 
 Let then $\zeta_{gen}$ denote the zeta-function of the generic $s$-monodromy $h_{gen}$ on $H_*(F_{s,t})$. Let $Z_i := B_i\cap F_{s,t}$ and let $\zeta_{Z_i}$ denote its zeta-function. Remark that $Z_i$ is the Milnor fiber of a transversal singularity to a non-isolated singularity, see \S \ref{ss:mon2} for explanations.
 
  For $n\ge2$, by using (\ref{eq:dual}), it follows that
  $\zeta_{Z_i}(\ft)$ is equal to the product $\zeta_{B_i\cap \bY_{s,t}}(\ft) \cdot \zeta_{B_i\cap \bY^\ity_{s,t}}^{-1}(\ft)$ of the zeta-functions of the boundary singularity $(\bY_{0,t}, \bY^\ity_{0,t})$ at the point $(p_i, 0, t)$. 
 
 \begin{proposition}\label{p:zeta}
\[ \zeta_{gen}(\ft) = (1-\ft)^{- \chi (F_{0,t})} \cdot \prod_{i} 
[\zeta_{B_i\cap \bY_{s,t}}(\ft) \cdot \zeta_{B_i\cap \bY^\ity_{s,t}}^{-1}(\ft)].\]
 \end{proposition} 
 \begin{proof}
 Note that the homology $\tilde H_*(F_{s,t})$ is concentrated in dimension 
$n-1$.
 The exact sequence of the pair $(F_{D,t}, F_{s,t})$
yields the following exact sequence on which the $s$-monodromy acts (for $n\ge 2$):
 \[ 0 \to H_{n}(F_{D,t}) \to H_{n}(F_{D,t}, F_{s,t}) \to 
H_{n-1}(F_{s,t}) \to H_{n-1}(F_{D,t}) \to 0.\]

By  (\ref{eq:sum}) and by Lemma \ref{l:ball}, we get $\prod_{i} \zeta_{Z_i}(\ft)= \zeta^{-1}_{(F_{D,t}, F_{s,t})}(\ft)$. As shown before, we also have $\zeta_{Z_i} = \zeta_{B_i\cap \bY_{s,t}} \cdot \zeta_{B_i\cap \bY^\ity_{s,t}}^{-1}$.
 Then apply the zeta-function to above exact sequence. 
 Since on $H_*(F_{D,t})$ the 
$s$-monodromy is trivial, and since $\chi(F_{D,t}) = \chi(F_{0,t})$, we get the claimed equality.
 \end{proof}

%%%%%%%%%%%%%%%%%%%%%%%%%%%%%%%

\subsection{The atypical case}

This case works exactly as the generic case and we get a similar formula for $\zeta_{aty}$, for a fixed value $t_0$:
\[\zeta_{aty}(\ft) = (1-\ft)^{- \chi (F_{0,t_0})} \cdot \prod_{i} [\zeta_{B_i\cap \bY_{s,t_0}}(\ft) \cdot \zeta_{B_i\cap \bY^\ity_{s,t_0}}^{-1}(\ft)].
\]
 Note that $F_{0,t_0}$ is now the atypical fiber and that 
the zeta-functions $\zeta_{B_i\cap \bY_{s,t_0}}$ and $\zeta_{B_i\cap \bY^\ity_{s,t_0}}$ are the ones of its corresponding boundary 
singularities.

%%%%%%%%%%%%%%%%%%%%%%%%%%%%%%%

\subsection{The global $s$-monodromy of $(F_{s,\bC}, F_{s,t})$.}

 The $s$-monodromy also acts on the homology exact sequence of the pair 
$(F_{s, \bC}, F_{s,t})$. Since $F_{s, \bC} = \bC^n$, the monodromy on it is trivial and we have 
  the isomorphism:  
  \[ H_*(F_{s, \bC}, F_{s,t}) \stackrel{\simeq}{\to} \tilde H_{*-1}(F_{s,t}), 
\]
   which is equivariant with respect to the $s$-monodromy action,
   concentrated in $*=n$.
   
 We may decompose the relative homology, by excision, into 
 a direct sum of submodules,  which correspond to the points where 
the discriminant $\Delta_P$ intersects the $t$-axis and at $t=\infty$.   We take within 
$\{s\} \times \bC$ a disjoint union of disks 
$\cup_j T_j$ connected by simple, non-intersecting paths to some general point 
$(s,t)$, with the following property: each disk $T_j$ contains all 
intersection points with $\Delta_P$ that tend, as $s\to 0$, to a single point 
of $\Delta_P \cap (\{ 0\} \times \bC)$ or tend to the infinity point $\{ 0\} \times \ity$ of 
the projective compactification $\{ 0\} \times \bP^1$. Similarly as in \cite{ST-m}, one may define a 
geometric monodromy on $F_{s, \bC}$, such that it acts on $F_{s, T_j}$, for any $j$, and on every fiber $F_{s, t}$ such that $t\not\in 
\cup_{j\ge 0} T_j$.
   By convention,  $T_0$ means the punctured disk centered at  $\{ 0\} \times 
\ity$; this corresponds of course to the intersections with $\Delta_{III}$.
   
  We get the following direct sum decomposition, compatible with the monodromy action in homology:
 \begin{equation}\label{eq:split}
  H_n(F_{s, \bC}, F_{s,t}) \stackrel{\simeq}{\to} \oplus_j H_n 
(F_{s,T_j},F_{s,t} ). 
\end{equation}

 In the source space we can further localise the vanishing homology to 
neigborhoods of the intersection points of $\Gamma_P$ with $\bY_{0,\bC}$.
The branches of $\Delta_P$ that meet $T_j$ can be projections of one or more branches of $\Gamma_P$, of different types. We 
assume from now on that each $T_j$ contains projections of a single type of branches, so that to each disk $T_j$ one may assign a well-defined type I, II or III.
 
The direct sum (\ref{eq:split}) can be decomposed into 3 terms, according to the type of $T_j$. Then the 
zeta-function splits accordingly too. We get:
 
%%%%%%%%%%%%%%%%%
 \begin{proposition}\label{p:zeta1}
 $\zeta_{gen}(\ft)= (1-\ft)^{-1+(-1)^n \mu(0)}  \zeta_{II}^{-1}(\ft) \cdot \zeta_{III}^{-1}(\ft)$.
 \end{proposition}  
\begin{proof} 
Comparing $\zeta_{gen}$ with the zeta-function of $h$ on 
 $H_n(F_{s, \bC}, F_{s,t})$ gives
 $\zeta_{gen}(\ft)= (1-\ft)^{-1}  \zeta_{I}^{-1}(\ft) \cdot \zeta_{II}^{-1}(\ft) \cdot 
\zeta_{III}^{-1}(\ft)$.

  For $\zeta_I$ we have the following formula:
  \[ \zeta_I (\ft) = (1-\ft)^{-\chi(F_{s,T_I},F_{s,t})} =  
(1-\ft)^{(-1)^{n+1}\mu(0))},\]
where $T_I$ denotes the union of type I disks and their corresponding paths. 
Indeed, 
for disks of type I, the $s$-monodromy on $H_n(F_{s,T_j})$ is the 
identity, and $\chi(F_{s,T_I},F_{s,t})=  (-1)^n \mu(0)$.
\end{proof}

%%%%%%%%%%%%%%%%%%%%%%%%%%%%%%%%%%%%
%%%%%%%%%%%%%%%%%%% 
\subsubsection{$s$-monodromy related to $\Delta_{II}$.}\label{ss:mon2}

 %%%%%%%%%%%%%%%%%%%%%%%%%%%%%%%%%%%%
For some point $(p,0,t_0)\in \Gamma_{II} \cap \bY^\ity_{0,T}$, let us take a Milnor ball $B$ and recall the notation $Y_{s,T}:= B \cap \bY_{s,T}$. This is the Milnor fiber 
of the germ at $(p,0,t_0)$ of the {\it non-isolated line singularity}\footnote{
Non-isolated line singularities and their monodromy have been studied by Dirk Siersma in several papers, see e.g. \cite{Si-top, Si-newton}.} $Y_{0,T}$ and $\zeta_{Y_{s,T}}$ will denote its monodromy. The reduced homology of $Y_{s,T}$ is concentrated in dimensions
$n-2$ and $n-1$. Slicing $Y_{0,T}$ by $\bY_{\bC,t}$ at some point $(p,0, t)$, for $t$ close to $t_0$, defines the {\em transversal singularity}.
Its Milnor fiber is $Z_t := B(p,0,t) \cap \bY_{s',t}$, where $B(p,0,t) \subset \bY$ is a small Milnor ball centered at $(p,0,t)$ and $s'\in \partial \bar D_t$, for some disc $D_t$ of radius smaller than that of $B_t$, such that $(D_t \times \{ t\} ) \cap \Gamma_P = \emptyset$.\footnote{
Remark that we used the notation $Z_i$ for the same object $Z_t$, in the beginning of \S \ref{s:zeta}.}
The reduced homology $\tilde H_*(Z_t, \bZ)$ is concentrated in dimension $n-2$. We denote by $\zeta_{Z_t}$ the zeta-function of the Milnor monodromy of the transversal singularity.

With these notations, we get the following formula for type II components: %%%%%%%%%%%%%%%%%%%%%%
\begin{proposition}\label{p:two}
 \[ \zeta_{II}(\ft) = (1-\ft)^{(-1)^{n-1}\lambda(0)} \prod [\zeta_{Y_{s,T}}(\ft) \cdot \zeta^{-1}_{Z_t}(\ft)], \] 
where the product is taken over all $T$'s of type II and, for each such $T$, over all points $\Gamma_{II} \cap \bY^\ity_{0,T}$.
 \end{proposition} 
%%%%%%%%%%%%%%%%%%%%%%
\begin{proof}
Note first that $Y_{s,t}:= B \cap \bY_{s,T}$ is the Milnor fiber
of the germ at $(p,0,t_0)$ of the hypersurface $\bY_{0,t_0}$. Then,  
for some fixed $T$ we have:
\[ \begin{array}{c}\zeta_{(F_{s,T},F_{s,t})} = \zeta_{F_{s,T}} 
\cdot \zeta^{-1}_{F_{s,t}} = 
 \zeta_{\bY_{s,T}} \cdot \zeta^{-1}_{\bY^{\ity}_{s,T}} \cdot 
\zeta^{-1}_{\bY_{s,t}} \cdot \zeta_{\bY^{\ity}_{s,t}} = \\
\prod [ \zeta_{Y_{s,T}} \cdot \zeta^{-1}_{Y_{s,t}} \cdot \zeta^{-1}_{Y^{\ity}_{s,T}}
\cdot \zeta_{Y^{\ity}_{s,t}}] 
 = \prod [\zeta_{Y_{s,T}} \cdot \zeta^{-1}_{Y_{s,t}}],
 \end{array}  \]
 where the product is taken over all points $\Gamma_{II} \cap \bY^\ity_{0,T}$.
 The cancellation in the last equality is due to 
the product structure $Y^{\ity}_{s,T} = Y^{\infty}_{s,t} \times T$.

  The zeta-function $\zeta_{Y_{s,t}}$ corresponds to the $s$-monodromy
  of the fiber $Y_{s,t}$ around the transversal singularity. This is not the Milnor monodromy of $Y_{s,t}$ since not all vanishing cycles of $Y_{s,t}$ vanish at the singularity $(p,0,t)$.
  In fact there is a geometric decomposition of the reduced homology 
into a direct sum  $\tilde H_*(Y_{s,t}, \bZ) = M\oplus N$, where 
$M$ is the module of cycles vanishing at  $(p,0,t)$ and $N$ denotes the module of cycles vanishing at the points of intersection $\Gamma_P \cap \bY_{t, \bC}$ which tend to $(p,0,t_0)$ as $t\to t_0$.
Note that $b_{n-1}(M) = \mu_{p, gen}(0)$ and that  $b_{n-1}(N) = \lambda_{p, t_0}(0)$, as shown in the proof of Theorem \ref{t:lambdaexchange} (see also Figure \ref{f:0}).

  Now we may apply a general Picard type theorem\footnote{a proof can be found in \cite[Prop. 2.2]{ST-mero}.} which in our case says that the action of the $s$-monodromy on some cycle in $N$ is the identity, modulo $M$. At the level of zeta-functions, this yields a product decomposition:
  \[   \zeta_{Y_{s,t}} (\ft) = (1-\ft)^{(-1)^n \lambda_{p,t_0}(0)}\cdot \zeta_{Z_t} (\ft).\]
 Taking the product over all points $\Gamma_{II} \cap \bY^\ity_{0,T}$
 and over all $T$, yields the claimed formula.  
\end{proof}
%%%%%%%%%%%%%%%%%%%

In terms of the corresponding Betti numbers, we can read the above proposition as follows:

\begin{equation}\label{eq:mu2}
 \mu_{II} = \sum_{(p,t_0)\in\Gamma_{II} \cap \bY^\ity_{0,T}} (b_n(Y_{0,T}, (p, t_0)) - b_{n-1}(Y_{0,T}, (p, t_0)) + \mu_{p,\gen}(0)) + \lambda(0), 
\end{equation}
where $b_n$ and $b_{n-1}$ are the Betti numbers of the Milnor fiber of the non-isolated line singularity $Y_{0,\bC}$ at the point $(p, 0, t_0)$.

%%%%%%%%%%%%%%%%%%%
\subsubsection{$s$-monodromy related to $\Delta_{III}$.}

%%%%%%%%%%%%%%%%%%%%%%%%%%%%%%%%

Since we can compute the zeta-function of the generic $s$-monodromy in two different ways, we are able to compute $\zeta_{III}$. Indeed, by comparing the formulas in Propositions  \ref{p:zeta} and \ref{p:zeta1}, we get:

\begin{equation}\label{c:zeta23}
 \zeta_{II}(\ft) \cdot  \zeta_{III}(\ft) = (1-\ft)^{(-1)^{n-1}\lambda (0)} \cdot 
\prod_{i} \zeta^{-1}_i(\ft).
\end{equation}

In terms of Milnor numbers, this yields:

\[ \mu_{II} + \mu_{III} = \lambda (0) + \sum \beta_i, \]

and so we find again the formula given by Lemma \ref{c:mu23}. 

From (\ref{c:zeta23}) and Proposition \ref{p:two}, we may deduce a formula for $\zeta_{III}$.

%%%%%%%%%%%%%%%%%%
%%%%%%%%%%%%%%%%%%%
 
\section{Examples}

%%%%%%%%%%%%%%%%%%%
%%%%%%%%%%%%%%%%%%%
We show by examples the various ways in which the $\mu \rightleftarrows \lambda$ exchange 
takes place during deformations $P$, in relation to the affine discriminant 
$\Delta_P$.
%%%%%%%%%%%%%%%%%%%

\begin{example}\label{e:broughton} 
Let $f = x^2y +x$, the Broughton example \cite{Bro} of a polynomial
 with no affine singularities and one singularity at infinity: $\mu = 0$ and $\lambda 
= 1$. The $\lambda$ corresponds to a jump of the local type
 of $\bX_t$ at infinity from $A_2$ to $A_3$, see \S \ref{s:inv}. Compare also to Remark \ref{r:2}.
  
Consider the deformation  $f_s = x^2y +x + sy^3$.
This is a GI deformation of a \fisi\  polynomial. 
 The topological bifurcation set $\Bif^{top}_P$ consists of the axis $\{ \sigma =0\}$ and 
the affine discriminant $\Delta_P$, which has equation $s= \alpha t^4$.
 
 We have that $\Delta_P = \Delta_{II}$ and $\Delta_{II}$ is tangent to
$\{ \sigma =0\}$. According to Theorem \ref{t:lambdaexchange}, this forces $\lambda(0) > 0$.
For $s\not=0$ we have: $\mu(s) = 4$ and $\lambda(s) = 0$, so, as $s\to 0$:
\[ \mu +\lambda = 4 + 0 \to 0 + 1.\]

\begin{figure}[hbtp]
\begin{center}
\epsfxsize=7cm
\leavevmode
\epsffile{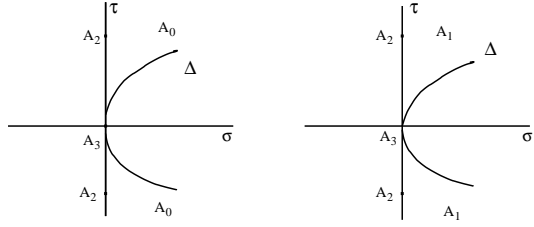}
\end{center}
\caption{{\em 
 $y^3$-deformation and $y$-deformation of Broughton's example}}
\label{f:1}   
\end{figure}

 In order to show the intermediate steps in this process, we consider the deformation, see Siersma-Smeltink paper \cite{SS}:
 \[  P(x,y; s_1, s_2, s_3) = x^2y + x + s_1 y + s_2 y^2 + s_3 y^3.\]
 Then the affine discriminant $\Delta_P$ has the following structure:
 \[ \begin{array}{lcl}
 s_3 \not= 0 & \ & \mu + \lambda = 4 + 0 \\
 s_3= 0,\  s_2 \not= 0 & \ & \mu + \lambda = 3 + 0 \\
 s_3= s_2 =0,\  s_1 \not= 0 & \ & \mu + \lambda = 2 + 0 \\
 s_3= s_2 = s_1= 0 & \ & \mu + \lambda = 0 + 1. 
 
\end{array}\]

In particular, the deformation $f_s = x^2y + x + s y$ is shown in Figure \ref{f:1}, right 
side. Here $\Delta_P = \Delta_{II}$ has equation $s= \alpha t^2$. 

For the zeta-function of the $s$-monodromy in the above $y^3$-deformation, Proposition \ref{p:zeta} gives the following formulas. Remark that the boundary singularity is of type $(A_2, A_1)$ at the point $[0:1]\in \bY^\ity_{0,t}$, for $t\not= 0$, and of type  $(A_3, A_1)$ when $t=0$.

\[ \zeta_{gen}(\ft) = \frac{1-\ft^3}{1-\ft}\cdot\frac{1-\ft^2}{1-\ft},\
\ \  \ \ \ \ \zeta_{aty}(\ft) = (1-\ft)^{-1}\cdot \frac{1-\ft^4}{1-\ft}\cdot\frac{1-\ft^2}{1-\ft}.\]

The formula for $\zeta_{II}$ in Proposition \ref{p:two}
needs knowledge on the non-isolated singularity of $\bY$ at the point $([0:1], 0, 0)$. According to de Jong's list \cite{dJ}, this is of type $F_1A_2$ and therefore its Milnor fiber is homotopy equivalent to a 2-sphere and its monodromy on $H_2$ is $-\id$. We get:

\[ \zeta_{II}(\ft) = (1-\ft)^{-1}\cdot \frac{1-\ft}{1-\ft^3}\cdot\frac{1}{1+\ft}.\]

This fits with (\ref{c:zeta23}) since $\zeta_{III}=1$ and the above result for $\zeta_{gen}$.
\end{example}
  
%%%%%%%%%%%%%%%%%%%%%%%%%%%%%%%%%%%%%%%%%%%%
\begin{example}\label{e:arnold} 
   Start with the polynomial $f = y^{k+1} + x^k$.  By homogenising and localising in the chart $x=1$, this corresponds to the fraction $A_k$ in  Arnol'd's list (Table \ref{t:1}), see \S \ref{s:arnold} for details.
   
    The hypersurface $\bX_t$ is non-singular and $k$-fold tangent to the hyperplane at infinity $H^\ity$.
   Let us study the following Yomdin deformation (in particular, a GI deformation): 
 \[ f_s = y^{k+1} + x^k + sx^{k+1}.\]
 Then $\bY_{s,t}$ intersects transversely $H^{\ity}$ for all $t$ and all $s\not= 0$. Moreover,
$\bY_{0,t}$ has no singularity at $H^{\ity}$ and is tangent to $H^{\ity}$ at the point $[1:0]\in H^{\ity}$.
 
  We get that $\Delta_I$ is the $s$-axis $\{ \tau = 0\}$, 
 $\Delta_{II}= \emptyset$ (and therefore $\lambda(0) = 0$) and that $\Delta_{III}$ has equation 
$s^k t= \alpha$ (see Figure \ref{f:2}).
  A slice $t= t_0$ will cut $\Gamma_{III}$
 at $k$ points, each of which representing a singularity of some $\bY_{s_i, t}$ of type $A_k$. 
 
 \begin{figure}[hbtp]
\begin{center}
\epsfxsize=3cm
\leavevmode
\epsffile{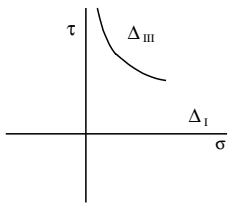}
\end{center}
\caption{{\em 
  Yomdin deformation of $A_k$.}}
\label{f:2}   
\end{figure}

\end{example}

Since at $H^\ity$ we have a boundary singularity of type $(A_0, A_k)$, the zeta-function of the generic $s$-monodromy is, after Proposition \ref{p:zeta}:
\[ \zeta_{gen}(\ft) = (1-\ft)^{-1-k+k^2} \cdot \frac{1-\ft^{k+1}}{1-\ft}.\]

%%%%%%%%%%%%%%%%%%%%%%%%%%%%%%%%%%%%%%%%%%%%
\begin{example}\label{e:mixed}(``Mixed example'').
\[  f_s = xy^4 + x^2y^2 + y + sx^5 .\]
 We get:
\[  \begin{array}{lcl}
 s= 0 & \ &  \mu + \lambda = 5+2 \\
 s\not= 0 &  \ & \mu + \lambda = 16 + 0.
 \end{array}\]
 When $s=0$, the value 2 of $\lambda(0)$ is got in connection to the jump
 $D_6 \to D_8$ at $t=0$. 
 
 The topological bifurcation set $\Bif_P^{top}$ consists of 3 branches: the axis $\sigma =0$ 
(since $\gamma(0) \not= \gamma(s)$,
 $\Delta_{I}$, $\Delta_{II}$ and $\Delta_{III}$. For  $\Delta_{II}$ we get the equation 
$s^2 = \alpha t^7$, and $\mu_{II} = 7$.
 For  $\Delta_{III}$ we get the equation $s^6 t^4= \beta$, and $\mu_{III} = 4$. The 
difference to the total $\gamma(s) =16$ is of course $\mu_I =5$, which corresponds to 
$\Delta_I$.

For the zeta-function of the $s$-monodromy, we apply Proposition \ref{p:zeta}, 
remarking that the boundary singularity at infinity is of type $(D_6, A_3)$, for generic $t$, and of type $(D_8, A_3)$ for $t=0$. We therefore get:

\[ \zeta_{gen}(\ft) = (1-\ft)^6\cdot \zeta_{D_6}(\ft)\cdot (1-\ft^4),\
\ \  \ \ \ \zeta_{aty}(\ft) = (1-\ft)^4\cdot\zeta_{D_8}(\ft)\cdot (1-\ft^4),\]
where $\zeta_{D_6}$ and $\zeta_{D_8}$ are the local zeta-functions of 
the $D_6$ and $D_8$ hypersurface germs (Arnol'd's classification \cite{Ar1}), respectively.
\end{example}
%%%%%%%%%%%%%%%%%%   
 %%%%%%%%%%%%%%%%%% 
\begin{example} \label{e:decrease}
In certain deformations, the invariant $\lambda$ may decrease. The following example is not a \fisi 
\ deformation since the initial polynomial $f_0$ is not \fisi \ (however, $f_0$ is isomorphic to a linear form): $f_s = x+ sx^2y +z^3$.
 We have: $ \lambda(s) = 2$, $\mu(s) =0$, for $s\not= 0$, and
   $ \lambda(0) = 0$, $\mu(0) =0$.
\end{example}   

%%%%%%%%%%%%%%%%%%
%%%%%%%%%%%%%%%%%
%%%%%%%%%%%%%

\section{On Arnol'd's classification}\label{s:arnold}
%%%%%%%%%%%%%%%%%
%%%%%%%%%%%%%
In his paper \cite{Ar2}, V.I. Arnol'd studied germs of fractions $P/Q$, where $P$ and $Q$ are polynomials, with respect to several equivalence relations such as the right-equivalence.
The relation to polynomial functions is as follows. Starting with a polynomial $f : \bC^n \to \bC$ of degree $d$ in coordinates $(x, y_2, \ldots , y_n)$,
we can introduce new affine coordinates $Y=\frac{y}{x}$, $Z=\frac{1}{x}$
in order to study points at infinity $[x:y:1] = [1:Y:Z]$.
For $f=t$ one gets $t = \frac{H(Y,Z)}{Z^d}$, where $H(Y,Z) = \tilde f(1, Y,Z)$ and $\tilde f(x, y, z)$ denotes the homogenised of $f(x,y)$ with the new variable $z$. 

According to Arnol'd's study, right-equivalence of fractions is very near  to right-equivalence of boundary singularities. In the above example, the boundary singularity is the one of the germ $H(Y,Z)$ with respect to the hyperplane $Z=0$.
  After  \cite{Ar2}, the list of simple fractions of type  $\frac{H(Y,Z)}{Z^d}$ is given in Table \ref{t:1} (where the letter $Q$ denotes a quadratic form in the rest of the $Y$-coordinates).
  
%%%%%%%%%%%%%%%%%%%%%%%%%%%%%%%%%%%%%%%
\begin{table}[hbtp]
\begin{center}
\epsfxsize=7cm
\leavevmode
\[
\begin{array}{|l|l|l|l|}
\hline
\text{type} & \text{function} & \text{boundary type} & \text{conditions}\\
\hline
A_0  & Y/Z^d & (A_0,A_k) &   \\
A_k  & (Y^{k+1} + Z + Q_2)/Z^d &(A_0,A_k)  & k \ge 1, \ d>1 \\
B_k  & (Y^{2} + Z^k + Q_2)/Z^d & (A_{k-1},A_1)  &   k \ge 2, \ d>k \\
C_k  &(Y^{k} + YZ + Q_2)/Z^d   & (A_1,A_{k-1})&   k \ge 3, \ d>1 \\
D_k  & (Y_1^{2}Y_2 + Y_2^{k-1} + Z + Q_3)/Z^d   &(A_0,D_k) &   k \ge 4, \ d>1 \\
E_6  & (Y_1^{3}+ Y_2^{4} + Z + Q_3)/Z^d & (A_0,E_6) &  d>1\\
E_7  & (Y_1^{3}+ Y_1Y_2^{3} + Z + Q_3)/Z^d & (A_0,E_7) &  d>1\\
E_8  & (Y_1^{3}+ Y_2^{5} + Z + Q_3)/Z^d & (A_0,E_8) &  d>1\\
F_4  & (Y^3 + Z^2 + Q_2)/Z^d & (A_2,A_2) &  d>2 \\
\hline
\end{array}
\]
\end{center}
\caption{{\em 
  Arnol'd's list of fractions of type $H/Z^d$ }}
\label{t:1}   
\end{table}
%%%%%%%%%%%%%%%%%%%%%%%%%%%%%%%%%%%%%%%
  
 All the simple germs in Arnol'd's 
list have $\lambda = 0$, so they are topologically trivial at infinity.  However, deforming 
these singularities yields an interesting behaviour. As we have shown before, we may consider for instance Yomdin deformations and their $s$-monodromy. Those are non-trivial. 

\begin{example}\label{e: final}
Let us treat the $A_k$ singularity, 
where $d= k+1$, namely $f = y^{k+1} + x^k$.
 Consider the following unfolding and its corresponding localisation:
 \[ \begin{array}{l}
 P(x,y;s_1, \ldots , s_k) = y^{k+1} + x^k + s_1x^2y^{k-1} + s_2x^3y^{k-2} + \cdots 
s_kx^{k+1}, \\
 y^{k+1} + x_0 + s_1y^{k-1} + s_2y^{k-2} + \cdots s_k = t x_0^{k+1}.
 \end{array}\]
 
 In the parameter space $\bC^k$ there exists a natural stratification of  the {\em bifurcation set} of $P$ (see Corollary \ref{l:bif}):
\[ 
 \Bif_P := \{ s\in \bC^k \mid b_{n-1}(G_s) < (d-1)^n \},
\]
corresponding to the number and multiplicity of the roots of the polynomial $y^{k+1} + 
s_1y^{k-1} + s_2y^{k-2} + \cdots s_k$. Each stratum of $\Bif_P$ yields a well-defined Betti 
number $b_1(G_s)$ of the general fiber $G_s$ of $f_s$.
The following 2-parameter deformation:
 \[  f_s = y^{k+1} + x^k + s_1x^{k+1} + s_2 x^{k}y\]
 
  yields a kind of Morsification, in the sense that,
 fixing $s_2$ close to $0$, we get $k$ values of $s_1$ for which 
 $f_s$ has simple tangencies at infinity.
  We send to Example \ref{e:arnold} for more details on the 1-parameter Yomdin deformation  $f_s = y^{k+1} + x^k + sx^{k+1}$.
 \end{example}
 %%%%%%%%%%%%%%%%%%
 We end by two remarks.
%%%%%%%%%%%%%%%
\begin{remark}\label{r:1}
 There is an important difference between looking to global polynomial functions and looking to the germ at some point at infinity.
  For example, take the $A_k$ fraction singularity from Arnol'd's list (Table \ref{t:1}). By changing coordinates, one gets $f = x^{d-k-1}y^{k+1} + x^{d-1}$, which can be viewed as a global polynomial function. When $d>k-1$, this has a non-isolated singularity, which is of course not in the chart $x=1$. However, one can get rid of non-isolated singularities  by finding more suitable representatives in the same right-equivalence class. For instance, in case of $A_1$:
\[  \frac{Y^2 + Z}{Z^d} \ \stackrel{\cR}{\sim} \ \frac{Y^2 + Z + Y^3Z}{Z^d},
\]
which corresponds, at the point $[1:0:0]$, to the equivalence:
\[  x^2y^2 + x^3 \sim x^2y^2 + x^3 +y^3.
\]
\end{remark}
%%%%%%%%%%%%%%%%%
\begin{remark}\label{r:2}
The polynomial $f= y + xy^2$ of degree $3$, which is Broughton's example (see Example \ref{e:broughton}), does not occur in Arnol'd's list (Table \ref{t:1}). Indeed, in coordinates at infinity we see the following equivalence:
\[  \frac{Y^2 + YZ^2}{Z^3} \ \stackrel{\cR}{\sim} \ \frac{Y^2 + Z^4}{Z^3},
\]
which is the $B_4$ fraction germ, but this is not simple, since here $d<k$. Another way to see the non-simplicity of $B_4$ is as follows. There is no inner modulus in case of $B_4$ in degree 3, but this exists for $B_3 :\  \frac{Y^2 + \alpha Z^3}{Z^3}$, since $\alpha$ is indeed an inner modulus. That $B_4$ is not simple follows now from adjacency. 
\end{remark} 
%%%%%%%%%%%%%%%%%%
%%%%%%%%%%%%%%% 

%

%%%%%%%%%%%%%%%%%%
%%%%%%%%%%%%%%%%%%
%%%%%%%%%%%%%%%%%%

\end{document}